\documentclass[12pt,english,leqno]{article}
\usepackage[T1]{fontenc}
\usepackage[latin1]{inputenc}
\usepackage{babel}

\makeatletter

\providecommand{\LyX}{L\kern-.1667em\lower.25em\hbox{Y}\kern-.125emX\@}

\usepackage{amsfonts,amsmath}
 \newtheorem{theorem}{Theorem}
 
 \newtheorem{lemma}[theorem]{Lemma}

 \makeatother
 
\begin{document}

 \title{Small time asymptotics in local limit theorems for Markov chains converging
to diffusions. \thanks{
 This research was supported by grant 436RUS113/845/0-1
from the Deutsche Forschungsgemeinschaft and by grants 04-01-00700
and 05-01-04004 from the Russian Foundation of Fundamental
Researches. The author worked on the paper during a visit at the
Laboratory of Probability Theory and Random Models of the
University Paris VI in 2005/2006. He is grateful for the
hospitality during his stay.
 }\\
 }
\author{Valentin KONAKOV \thanks{
 Laboratory of Probability and Statistics, Central Economics Mathematical Institute, Academy of
 Sciences,
 Nahimovskii av. 47, 117418 Moscow, Russia, E mail:
  kv24@mail.ru}}

 \date{\today{}}

 \maketitle
 \noindent We consider triangular arrays of Markov chains that converge weakly  to a
diffusion process. Local limit theorems for transition densities
are proved. The observation time $[0,T]$ \ may be fixed or  $\lim
$ $_{n\rightarrow \infty }T=0,$ where $\ nh=T$ \ and $\ h$ \ is a
mesh between two neighboring observation points. \vskip .2in
\textsl{1991 MSC:} primary
 62G07, secondary 60G60


 \newpage

{\protect\bigskip \bf 1. Introduction and main results}

Let $n\geq2$ and $h>0$ be such that $T=nh\leq1.$ Suppose that $%
q\left(t,x,\cdot\right),$ $\left(t,x\right)\in\lbrack0,1]\times\Bbb{R}^{d}$
is a given family of densities on $\Bbb{R}^{d}$ and $m$ is a function from $%
[0,1]\times\Bbb{R}^{d}$ into $\Bbb{R}^{d}.$ We shall impose the following
conditions

\begin{description}
\item[A1]  $\int_{\Bbb{R}^{d}}yq\left( t,x,y\right) dy=0,\;0\leq t\leq 1,\
x\in \Bbb{R}^{d}.$

\item[A2]  There exists positive constants $\sigma _{\ast }$ and $\sigma
^{\ast }$ such that the covariance matrix $\sigma \left( t,x\right) =\int_{%
\Bbb{R}^{d}}yy^{T}q\left( t,x,y\right) dy$ satisfies
\[
\sigma _{\ast }\leq \theta ^{T}\sigma \left( t,x\right) \theta \leq \sigma
^{\ast },
\]
for all $\left\| \theta \right\| =1$ and $t\in \lbrack 0,1]$ $x\in \Bbb{R}%
^{d}.$

\item[A3]  There exists a positive integer $S^{\prime }$ and a real
nonnegative function $\psi \left( y\right) ,$ $y\in \Bbb{R}^{d}$ satisfying $%
\sup_{y\in \Bbb{R}^{d}}\psi \left( y\right) <\infty $ and $\int_{\Bbb{R}%
^{d}}\left\| y\right\| ^{S}\psi \left( y\right) dy<\infty ,$ with $%
S=2dS^{\prime }+4,$ such that
\[
\left| D_{y}^{\nu }q\left( t,x,y\right) \right| \leq \psi \left( y\right)
,\;t\in \lbrack 0,1],\;x,y\in \Bbb{R}^{d}\;\left| \nu \right| =0,1,2,3,4
\]
and
\[
\left| D_{x}^{\nu }q\left( t,x,y\right) \right| \leq \psi \left( y\right)
,\;t\in \lbrack 0,1],\;x,y\in \Bbb{R}^{d}\;\left| \nu \right| =0,1,2.
\]
Furthemore, for all $x\in \Bbb{R}^{d}$ it holds
\[
\int_{\Bbb{R}^{d}}\left| q\left( t,x,y\right) -q\left( t^{\prime
},x,y\right) \right| dy\rightarrow 0,
\]
as $\left| t-t^{\prime }\right| \rightarrow 0.$

\item[B1]  The functions $m\left( t,x\right) $ and $\sigma \left( t,x\right)
$ and their first and second derivatives w.r.t. $t$ and $x$ are continuous
and bounded uniformly in $t$ and $x.$ All these functions are Lipschitz
continuous with respect to $x$ with a Lipschitz constant that does not
depend on $t.$ Furthemore, $D^{\nu }\sigma \left( t,x\right) $ exists for $%
\left| \nu \right| =6$ \ and \ is Holder continuous w.r.t. $x$ with positive
exponent and a constant that does not depend on $t.$
\end{description}

Consider a family of Markov processes in $\Bbb{R}^{d}$ of the following form
\begin{equation}
X_{k+1,h}=X_{k,h}+m\left( kh,X_{k,h}\right) h+\sqrt{h}\xi
_{k+1,h},\;X_{0,h}=x\in \Bbb{R}^{d},\;k=0,...,n-1,  \label{I-001}
\end{equation}
where $\left( \xi _{i,h}\right) _{i=1,...,n}$ is an innovation sequence
satisfying the Markov assumption:\ the conditional distribution of $\xi
_{k+1,h}$ given $X_{k,h}=x_{k},...,X_{0,h}=x_{0}$ depends only on $%
X_{k,h}=x_{k}$ and has conditional density $q\left( kh,x_{k},\cdot \right) .$
The conditional covariance matrix corresponding to this density is $\sigma
(kh,x_{k}).$ The transition densities of $\ \left( X_{i,h}\right)
_{i=1,...,n}$ are denoted by $p_{h}\left( 0,kh,x,\cdot \right) .$

We shall consider the process (\ref{I-001}) as an approximation to the
following stochastic differential equation in $\Bbb{R}^{d}:$%
\[
dY_{s}=m\left( s,Y_{s}\right) ds+\Lambda \left( s,Y_{s}\right)
dW_{s},\;Y_{0}=x\in \Bbb{R}^{d},\;s\in \lbrack 0,T],
\]
where $\left( W_{s}\right) _{s\geq 0}$ \ is the standard Wiener process and $%
\Lambda $ is a symmetic positive definite $d\times d$ \ matrix such that $%
\Lambda \left( s,y\right) \Lambda \left( s,y\right) ^{T}=\sigma \left(
s,y\right) .$ The conditional density of $Y_{t},$ given $Y_{0}=x$ is denoted
by $p\left( 0,t,x,\cdot \right) .\,\,$

Konakov and Mammen (2000) obtained a nonuniform rate of convergence for the
difference $p_{h}\left( 0,T,x,\cdot \right) -p\left( 0,T,x,\cdot \right) $
as $n\rightarrow \infty $ in the case $T\asymp 1.$ It is the goal of the
present note to get an analogous result in the case $T=o\left( 1\right) .$
The following theorem is our main result.

\begin{theorem}
Let $h>0,$ $n\geq 2,$ $T=nh.$ Assume (A1-A3) and (B1). Then, as $%
n\rightarrow \infty $ and $T\rightarrow 0,$
\[
\sup_{x,y}Q_{\sqrt{T}}\left( y-x\right) \left| p_{h}\left( 0,T,x,y\right)
-p\left( 0,T,x,y\right) \right| =O\left( n^{-1/2}\right) ,
\]
where the constant in $O\left( \cdot \right) $ does not depend on $h$ and
\[
Q_{\delta }\left( u\right) =\delta ^{d}\left( 1+\left\| \frac{u}{\delta }%
\right\| ^{2S^{\prime }-2}\right) .
\]
\end{theorem}

{\protect\bigskip \bf 2. Parametrix method for diffusions}

For any $s\in (0,T),$ $x,y\in \Bbb{R}^{d}$ we consider an additional family
of \ ''frozen'' diffusion processes
\[
d\widetilde{Y}_{t}=m\left( t,y\right) dt+\Lambda \left( t,y\right) dW_{t},\;%
\widetilde{Y}_{s}=x,\;s\leq t\leq T.
\]
Let $\widetilde{p}^{y}\left( s,t,x,\cdot \right) $ be the conditional
density of $\widetilde{Y}_{t},$ given $\widetilde{Y}_{s}=x.$ In the sequel
for any $z$ we shall denote \ $\widetilde{p}\left( s,t,x,z\right) =%
\widetilde{p}^{z}\left( s,t,x,z\right) ,$ where the variable $z$ acts here
twice: as the arument of the density and as defining quantity of the process
$\widetilde{Y}_{t}.$

The transition densities $\widetilde{p}$ can be computed explicitly

\begin{eqnarray*}
\widetilde{p}\left( s,t,x,y\right) &=&\left( 2\pi \right) ^{-d/2}\left( \det
\sigma \left( s,t,y\right) \right) ^{-1/2} \\
&&\times \exp \left( -\frac{1}{2}\left( y-x-m\left( s,t,y\right) \right)
^{T}\sigma ^{-1}\left( s,t,y\right) \left( y-x-m\left( s,t,y\right) \right)
\right) ,
\end{eqnarray*}
where
\[
\sigma \left( s,t,y\right) =\int_{s}^{t}\sigma \left( u,y\right)
du,\;\;\;\;\;m\left( s,t,y\right) =\int_{s}^{t}m\left( u,y\right) du.
\]
We need the following kernel
\begin{eqnarray*}
H\left( s,t,x,y\right) &=&\frac{1}{2}\sum_{i,j=1}^{d}\left( \sigma
_{ij}\left( s,x\right) -\sigma _{ij}\left( s,y\right) \right) \frac{\partial
^{2}\widetilde{p}\left( s,t,x,y\right) }{\partial x_{i}\partial x_{j}} \\
&&+\sum_{i=1}^{d}\left( m_{i}\left( s,x\right) -m_{i}\left( s,y\right)
\right) \frac{\partial \widetilde{p}\left( s,t,x,y\right) }{\partial x_{i}}.
\end{eqnarray*}
Introduce the convolution type binary operation $\otimes :$%
\[
f\otimes g\left( s,t,x,y\right) =\int_{s}^{t}du\int_{\Bbb{R}^{d}}f\left(
s,u,x,z\right) g\left( u,t,z,y\right) dz.
\]
$k$-fold convolution of $H$ is denoted by $H^{\left( k\right) }.$ The
following results are taken from Konakov and Mammen (2000).\

\begin{lemma}
\label{Lemma-KM-1} Let $0\leq s<t\leq T.$ It holds
\[
p\left( s,t,x,y\right) =\sum_{r=0}^{\infty }\widetilde{p}\otimes H^{\left(
r\right) }\left( s,t,x,y\right) .
\]
\end{lemma}

\begin{lemma}
\label{Lemma-KM-2} Let $0\leq s<t\leq T.$ There exist constants $C$ \ and $%
C_{1}$ such that
\[
\left| H\left( s,t,x,y\right) \right| \leq C_{1}\rho ^{-1}\phi _{C,\rho
}\left( y-x\right)
\]
and
\[
\left| \widetilde{p}\otimes H^{\left( r\right) }\left( s,t,x,y\right)
\right| \leq C_{1}^{r+1}\frac{\rho ^{r}}{\Gamma \left( 1+\frac{r}{2}\right) }%
\phi _{C,\rho }\left( y-x\right) ,
\]
where $\rho ^{2}=t-s,$ $\phi _{C,\rho }\left( u\right) =\rho ^{-d}\phi
_{C}\left( u/\rho \right) $ and
\[
\phi _{C}\left( u\right) =\frac{\exp \left( -C\left\| u\right\| ^{2}\right)
}{\int \exp \left( -C\left\| v\right\| ^{2}\right) dv}.
\]
\end{lemma}

{\protect\bigskip \bf 3. Parametrix method for Markov chains}

For any $0\leq jh\leq T,$ $x,y\in \Bbb{R}^{d}$ we consider an additional
family of \ ''frozen'' \ Markov chains defined for $jh\leq ih\leq T$ \ as
\begin{equation}
\widetilde{X}_{i+1,h}=\widetilde{X}_{i,h}+m\left( ih,y\right) h+\sqrt{h}%
\widetilde{\xi }_{i+1,h},\;\widetilde{X}_{j,h}=x\in \Bbb{R}^{d},\,\,j\leq
i\leq n,  \label{Froz001}
\end{equation}
where $\widetilde{\xi }_{j+1,h},...,\widetilde{\xi }_{n,h}$ is an innovation
sequence such that the conditional density of $\widetilde{\xi }_{i+1,h}$
given $\widetilde{X}_{i,h}=x_{i},...,\widetilde{X}_{0,h}=x_{0}$ equals to $%
q\left( ih,y,\cdot \right) .$ Let us introduce the infinitesimal operators
corresponding to Markov chains (\ref{I-001}) and (\ref{Froz001})
respectively,
\begin{eqnarray*}
&&L_{h}f\left( jh,kh,x,y\right)\\
&& =h^{-1}\left( \int p_{h}\left( jh,\left( j+1\right)
h,x,z\right) f\left( \left( j+1\right) h,kh,z,y\right) dz-f\left(
\left( j+1\right) h,kh,z,y\right) \right)
\end{eqnarray*}
and
\begin{eqnarray*}
&&\widetilde{L}_{h}f\left( jh,kh,x,y\right)\\
&& =h^{-1}\left( \int \widetilde{p}%
_{h}^{y}\left( jh,\left( j+1\right) h,x,z\right) f\left( \left( j+1\right)
h,kh,z,y\right) dz-f\left( \left( j+1\right) h,kh,z,y\right) \right) ,
\end{eqnarray*}
where $\widetilde{p}_{h}^{y}\left( jh,j^{\prime }h,x,\cdot \right)
$ denotes
the conditional density of $\widetilde{X}_{j^{\prime },h}$ given $\widetilde{%
X}_{j,h}=x.$ As before for any $z$ denote $\widetilde{p}_{h}\left(
jh,j^{\prime }h,x,z\right) =\widetilde{p}_{h}^{z}\left( jh,j^{\prime
}h,x,z\right) ,$ where the variable $z$ acts here twise: as the arument of
the density and as defining quantity of the process $\widetilde{X}_{i,h}.$
For technical convenience the terms $f\left( \left( j+1\right)
h,kh,z,y\right) $ on the right hand side of $\ L_{h}f$ \ and $\widetilde{L}%
_{h}f$ \ \ appear instead of $f\left( jh,kh,z,y\right) .$

In analogy with the definition of $H$ we put, for $k>j,$%
\[
H_{h}\left( jh,kh,x,y\right) =\left( L_{h}-\widetilde{L}_{h}\right)
\widetilde{p}_{h}\left( jh,kh,x,y\right) .
\]
We also shall use the convolution type binary operation $\otimes _{h}:$%
\[
g\otimes _{h}f\left( jh,kh,x,y\right) =\sum_{i=j}^{k-1}h\int_{\Bbb{R}%
^{d}}g\left( jh,ih,x,z\right) f\left( ih,kh,z,y\right) dz,
\]
where $0\leq j<k\leq n.$ Write $g\otimes _{h}H_{h}^{\left( 0\right) }=g$ \
and $g\otimes _{h}H_{h}^{\left( r\right) }=\left( g\otimes _{h}H_{h}^{\left(
r-1\right) }\right) \otimes _{h}H_{h}$ \ for $r=1,...,n.$ For the higher
order convolutions we use the convention $\sum_{i=j}^{l}=0$ for $l<j.$ One
can show the following analog of the ''parametrix'' \ \ expansion for $p_{h}$
{[}see Konakov and Mammen (2000){]}.

\begin{lemma}
\label{Lemma-KM-3} Let $0\leq jh<kh\leq T.$ It holds
\[
p_{h}\left( jh,kh,x,y\right) =\sum_{r=0}^{k-j}\widetilde{p}_{h}\otimes
_{h}H_{h}^{\left( r\right) }\left( jh,kh,x,y\right) ,
\]
where
\[
p_{h}\left( jh,jh,x,y\right) =\widetilde{p}_{h}\left( kh,kh,x,y\right)
=\delta \left( y-x\right)
\]
and $\delta $ is the Dirac delta symbol.
\end{lemma}

{\protect\bigskip \bf 4. Auxiliary statements}

In this section we collect several bounds to be used later in the proofs.
Throughout the proofs $C,C_{1},C_{2},...$ denote generic constants possibly
different in different places.

We give some tools that are useful for comparison of the expansions of $p$
and $p_{h}$ (see Lemmas \ref{Lemma-KM-1}) and \ref{Lemma-KM-3}
respectively). Since $p$ and $p_{h}$ are written in terms of $\widetilde{p}%
^{y}$ and $\widetilde{p}_{h}^{y}$ it is essential first to bound the
difference $\widetilde{p}_{h}^{y}-\widetilde{p}^{y}.$ For this we use some
bounds from Konakov and Mammen (2000) which are proved for $T=1,$ but remain
true also for $T\rightarrow0.$ The following lemmas are slight modifications
of Lemmas 3.8, 3.9, 3.11, 3.12, 3.13 from Konakov and Mammen (2000) and
therefore the proofs will be not detailes here.

To simplify the notation let $\rho=\sqrt{\left(k-j\right)h}.$ Denote $%
\zeta_{\rho}\left(z\right)=\rho^{-p}\zeta\left(z/\rho\right),$ where
\[
\zeta\left(z\right)=\frac{\left(1+\left\Vert z\right\Vert ^{S-4}\right)^{-1}%
}{\int\left(1+\left\Vert z^{\prime}\right\Vert ^{S-4}\right)^{-1}dz^{\prime}}%
.
\]

\begin{lemma}
\label{Lemma-KM-4}Assume Conditions (A1-A3) and (B1). Then for $0\leq
j<k\leq n$ and all $x,y\in \Bbb{R}^{d}$ it holds
\[
\left| \widetilde{p}_{h}\left( jh,kh,x,y\right) -\widetilde{p}\left(
jh,kh,x,y\right) \right| \leq Ch^{1/2}\rho ^{-1}\zeta _{\rho }\left(
y-x\right) .
\]
\end{lemma}

For $j=0,...,k-2,$ let
\[
K_{h}\left(jh,kh,x,y\right)=\left(L-\widetilde{L}\right)\widetilde{p}%
_{h}\left(jh,kh,x,y\right)
\]
and
\begin{eqnarray*}
M_{h}\left(jh,kh,x,y\right) & = &
3h^{1/2}\sum_{\left|\nu\right|=3}\sum_{\left|\mu\right|=1}\int_{\Bbb{R}%
^{d}}d\theta\int_{0}^{1}d\delta
D_{y}^{\mu}q\left(jh,y,\theta\right)\left(x-y\right)^{\mu} \\
& & \times\frac{\theta^{\nu}}{\nu!}D_{x}^{\nu}\widetilde{p}%
_{h}\left(\left(j+1\right)h,kh,x+\delta\theta
h^{1/2},y\right)\left(1-\delta\right)^{2}.
\end{eqnarray*}
If $j=k-1$ define
\[
K_{h}\left(\left(k-1\right)h,kh,x,y\right)=0\,\,\,\,\,\,\mathrm{and\,\,}%
\,\,\,\,\ M_{h}\left(\left(k-1\right)h,kh,x,y\right)=0.
\]

\begin{lemma}
\label{Lemma-KM-5} Assume Conditions (A1-A3) and (B1). Then for $0\leq
j<k\leq n$ and all $x,y\in \Bbb{R}^{d}$ it holds
\[
\left| H_{h}\left( jh,kh,x,y\right) -K_{h}\left( jh,kh,x,y\right)
-M_{h}\left( jh,kh,x,y\right) \right| \leq Ch^{1/2}\rho ^{-1}\zeta _{\rho
}\left( y-x\right) ,
\]
where $\zeta _{p}$ is defined above.
\end{lemma}

Denote $\xi_{\rho}\left(z\right)=\rho^{-p}\xi\left(z/\rho\right),$ where
\[
\xi\left(z\right)=\frac{\left(1+\left\Vert z\right\Vert
^{2S^{\prime}-2}\right)^{-1}}{\int\left(1+\left\Vert z^{\prime}\right\Vert
^{2S^{\prime}-2}\right)^{-1}dz^{\prime}}.
\]

\begin{lemma}
\label{Lemma-KM-6}Assume Conditions (A1-A3) and (B1). Then for $r=1,2,...,$ $%
0\leq j<k\leq n$ and all $x,y\in \Bbb{R}^{d}$ it holds
\[
\left| \widetilde{p}_{h}\otimes _{h}H_{h}^{\left( r\right) }\left(
jh,kh,x,y\right) \right| \leq \frac{C^{r+1}\rho ^{r}}{\Gamma \left( 1+\frac{r%
}{2}\right) }\xi _{\rho }\left( x-y\right) .
\]
\end{lemma}

\begin{lemma}
\label{Lemma-KM-7}Assume Conditions (A1-A3) and (B1). Then for $0\leq
j<k\leq n$ and all $x,y\in \Bbb{R}^{d}$ it holds
\[
p_{h}\left( jh,kh,x,y\right) =\sum_{r=0}^{k-j}\left( \widetilde{p}%
_{h}\otimes _{h}\left( M_{h}+K_{h}\right) ^{\left( r\right) }\right) \left(
jh,kh,x,y\right) +R_{1}\left( jh,kh,x,y\right) ,
\]
where
\[
\left| R_{1}\left( jh,kh,x,y\right) \right| \leq Ch^{1/2}\rho \xi _{\rho
}\left( y-x\right) .
\]
\end{lemma}

\begin{lemma}
\label{Lemma-KM-8}Assume Conditions (A1-A3) and (B1). Then for $0\leq
j<k\leq n$ and all $x,y\in \Bbb{R}^{d}$ it holds
\[
p_{h}\left( jh,kh,x,y\right) =\sum_{r=0}^{k-j}\left( \widetilde{p}\otimes
_{h}\left( M_{h}+K_{h}\right) ^{\left( r\right) }\right) \left(
jh,kh,x,y\right) +R_{2}\left( jh,kh,x,y\right) ,
\]
where
\[
\left| R_{2}\left( jh,kh,x,y\right) \right| \leq Ch^{1/2}\rho ^{-1}\xi
_{\rho }\left( y-x\right) .
\]
\end{lemma}

{5. Proof of the main result}

From Lemmas \ref{Lemma-KM-1} and \ref{Lemma-KM-2} we get as $n\rightarrow
\infty $ and $T=nh\rightarrow 0$%
\begin{equation}
p\left( 0,T,x,y\right) =\sum_{r=0}^{n}\left( \widetilde{p}\otimes H^{\left(
r\right) }\right) \left( 0,T,x,y\right) +T^{-d/2}\exp \left( -\frac{C\left\|
y-x\right\| ^{2}}{T}\right) o\left( T^{n}e^{-n}\right) .  \label{MR-001}
\end{equation}
Furthemore, Lemma \ref{Lemma-KM-8} implies that
\begin{equation}
p_{h}\left( 0,T,x,y\right) =\sum_{r=0}^{n}\left( \widetilde{p}\otimes
_{h}\left( M_{h}+K_{h}\right) ^{\left( r\right) }\right) \left(
0,T,x,y\right) +\xi _{\sqrt{T}}\left( y-x\right) O\left( n^{-1/2}\right) .
\label{MR-002}
\end{equation}
Because of (\ref{MR-001}) and (\ref{MR-002}) for the statement of the
theorem it remains to show that
\begin{eqnarray}
&&\left| \sum_{r=0}^{n}\left( \widetilde{p}\otimes H^{\left( r\right)
}\right) \left( 0,T,x,y\right) -\sum_{r=0}^{n}\left( \widetilde{p}\otimes
_{h}\left( M_{h}+K_{h}\right) ^{\left( r\right) }\right) \left(
0,T,x,y\right) \right|  \nonumber \\
&&\;\;\;\;\;\;\;\;\;\;\;\;\;\;\;\;\;\;\;\;\;\;\;\;\;\;\;\;\;\;\;\;\;\;\;\;\;%
\;=\xi _{\sqrt{T}}\left( y-x\right) O\left( n^{-1/2}\right) .  \label{MR-003}
\end{eqnarray}
For the proof of (\ref{MR-003}) note that
\begin{eqnarray}
&&\left| \sum_{r=0}^{n}\left( \widetilde{p}\otimes H^{\left( r\right)
}\right) \left( 0,T,x,y\right) -\sum_{r=0}^{n}\left( \widetilde{p}\otimes
_{h}\left( M_{h}+K_{h}\right) ^{\left( r\right) }\right) \left(
0,T,x,y\right) \right| \;  \nonumber \\
\;\;\;\; &\leq &S_{1}+S_{2}+S_{3},  \label{MR-004}
\end{eqnarray}
where
\begin{eqnarray*}
S_{1} &=&\left| \sum_{r=0}^{n}\left( \widetilde{p}\otimes H^{\left( r\right)
}\right) \left( 0,T,x,y\right) -\sum_{r=0}^{n}\left( \widetilde{p}\otimes
_{h}H^{\left( r\right) }\right) \left( 0,T,x,y\right) \right| , \\
S_{2} &=&\left| \sum_{r=0}^{n}\left( \widetilde{p}\otimes _{h}H^{\left(
r\right) }\right) \left( 0,T,x,y\right) -\sum_{r=0}^{n}\left( \widetilde{p}%
\otimes _{h}\left( M_{h}+H\right) ^{\left( r\right) }\right) \left(
0,T,x,y\right) \right| , \\
S_{3} &=&\left| \sum_{r=0}^{n}\left( \widetilde{p}\otimes _{h}\left(
M_{h}+H\right) ^{\left( r\right) }\right) \left( 0,T,x,y\right)
-\sum_{r=0}^{n}\left( \widetilde{p}\otimes _{h}\left( M_{h}+K_{h}\right)
^{\left( r\right) }\right) \left( 0,T,x,y\right) \right| .
\end{eqnarray*}
For $S_{1},$ $S_{2},$ $S_{3}$ we will show the following estimates
\begin{equation}
S_{k}=Q_{\sqrt{T}}\left( y-x\right) O\left( h^{1/2}\right) ,\;k=1,2,3.
\label{MR-005}
\end{equation}
We shall prove (\ref{MR-005}) for $k=1.$ The term $S_{1}$ corresponds to the
passage from the continuous time to the lattice time. The fact that $%
T\rightarrow 0$ implies that integrands involved in the convolutions $%
\otimes $ and $\otimes _{h}$ become asymptotically degenerate and therefore
more accurate estimates are required than those in Konakov and Mammen
(2000). We will develop here the details for these bounds.

We start from the recurrence relations for $r=1,2,3,...$%
\[
\left(\widetilde{p}\otimes
H^{\left(r\right)}\right)\left(0,jh,x,y\right)-\left(\widetilde{p}%
\otimes_{h}H^{\left(r\right)}\right)\left(0,jh,x,y\right)\qquad\;\;\;\;
\]

\[
=\left[ \left( \widetilde{p}\otimes H^{\left( r-1\right) }\right) \otimes
H-\left( \widetilde{p}\otimes H^{\left( r-1\right) }\right) \otimes _{h}H%
\right] \left( 0,jh,x,y\right) \qquad
\]
\begin{equation}
\,\,\,\,+\left[ \left( \widetilde{p}\otimes H^{\left( r-1\right) }\right) \
-\left( \widetilde{p}\otimes _{h}H^{\left( r-1\right) }\right) \right]
\otimes _{h}H\left( 0,jh,x,y\right) .\qquad  \label{MR-006}
\end{equation}
By summing up the identities in (\ref{MR-006}) from $r=1$ to $\infty $ and
by using the linearity of the operations $\otimes $ and $\otimes _{h}$ we
get
\[
(p-p^{d})\left( 0,jh,x,y\right) =\left( p\otimes H-p\otimes _{h}H\right)
\left( 0,jh,x,y\right)
\]
\begin{equation}
+(p-p^{d})\otimes _{h}H\left( 0,jh,x,y\right)  \label{MR-007}
\end{equation}
where we put
\begin{equation}
p^{d}(ih,i^{\prime }h,x,y)=\sum_{r=0}^{\infty }(\widetilde{p}\otimes
_{h}H^{(r)})(ih,i^{\prime }h,x,y).  \label{MR-008}
\end{equation}
By iterative application of (\ref{MR-007}) we obtain
\[
(p-p^{d})\left( 0,jh,x,y\right) =\left( p\otimes H-p\otimes _{h}H\right)
\left( 0,jh,x,y\right)
\]
\begin{equation}
+\left( p\otimes H-p\otimes _{h}H\right) \otimes _{h}\Phi \left(
0,jh,x,y\right) ,\qquad \,\,\,\,\qquad \;\;  \label{MR-009}
\end{equation}
where $\Phi (ih,i^{\prime }h,z,z^{\prime })=H(ih,i^{\prime }h,z,z^{\prime
})+H\otimes _{h}H(ih,i^{\prime }h,z,z^{\prime })+...=\sum_{r=1}^{\infty
}H^{(r)}(ih,i^{\prime }h,z,z^{\prime }).$

By the Taylor expansion we have
\[
\left(p\otimes
H-p\otimes_{h}H\right)(0,jh,x,z)\qquad\qquad\qquad\qquad\qquad\qquad\;\;
\]
\[
=\sum_{i=0}^{j-1}\int_{ih}^{(i+1)h}du\int_{R^{d}}\left[\lambda\left(u%
\right)-\lambda\left(ih\right)\right]dv\qquad\qquad\qquad\qquad\;\;\,\,
\]
\[
=\sum_{i=0}^{j-1}\int_{ih}^{(i+1)h}(u-ih)du\int_{R^{d}}\lambda^{%
\prime}(ih)dv\qquad\qquad\qquad\qquad\qquad
\]
\begin{equation}
+\frac{1}{2}\sum_{i=0}^{j-1}\int_{ih}^{(i+1)h}(u-ih)^{2}\int_{0}^{1}(1-%
\delta)\int_{R^{d}}\lambda^{\prime\prime}(s)\mid_{s=s_{i}}dvd\delta du,
\label{MR-010}
\end{equation}
where $\lambda\left(u\right)=p(0,u,x,v)H(u,jh,v,z),$ $s_{i}=s_{i}(u,i,%
\delta,h)=ih+\delta(u-ih).$

Note that
\[
\int_{R^{d}}\lambda^{\prime}(ih)dv=\int_{R^{d}}\frac{\partial}{\partial s}%
p(0,s,x,v)\mid_{s=ih}H(ih,jh,v,z)dv\qquad\;\;\;\;\;
\]
\[
+\int_{R^{d}}p(0,ih,x,v)\frac{\partial}{\partial s}H(s,jh,v,z)\mid_{s=ih}dv=%
\int_{R^{d}}L^{t}p(0,ih,x,v)\;\;
\]
\[
\times(L-\widetilde{L})\widetilde{p}(ih,jh,v,z)dv-\int_{R^{d}}p(0,ih,x,v)[(L-%
\widetilde{L})\widetilde{L}\widetilde{p}(ih,jh,v,z)
\]
\[
-H_{1}(ih,jh,v,z)]dv=\int_{R^{d}}p(0,ih,x,v)H_{1}(ih,jh,v,z)dv\qquad\qquad\;%
\;
\]
\begin{equation}
+\int_{R^{d}}p(0,ih,x,v)(L^{2}-2L\widetilde{L}+\widetilde{L}^{2})\widetilde{p%
}(ih,jh,v,z)dv,\qquad\qquad\qquad\;\;  \label{MR-011}
\end{equation}
where $H_{1}(s,t,v,z)$ is defined below in (\ref{MR-015}). We get from (\ref
{MR-011})
\[
\sum_{i=0}^{j-1}\int_{ih}^{(i+1)h}(u-ih)du\int_{R^{d}}\lambda^{\prime}(ih)dv=%
\frac{h}{2}(p\otimes_{h}H_{1})(0,jh,x,z)
\]
\begin{equation}
+\frac{h}{2}(p\otimes_{h}A_{0})(0,jh,x,z),\qquad\qquad\qquad\qquad\qquad%
\qquad\qquad\;\;\;  \label{MR-012}
\end{equation}
where $A_{0}(s,jh,v,z)=(L^{2}-2L\widetilde{L}+\widetilde{L}^{2})\widetilde{p}%
(s,jh,v,z).$ The direct calculation shows that
\[
A_{0}(s,jh,v,z)=\frac{1}{4}\sum_{p,q,r,l=1}^{d}(\sigma_{pq}(s,v)-%
\sigma_{pq}(s,z))(\sigma_{rl}(s,v)-\sigma_{rl}(s,z))\qquad\;
\]
\[
\times\frac{\partial^{4}\widetilde{p}(s,jh,v,z)}{\partial v_{p}\partial
v_{q}\partial v_{r}\partial v_{l}}+\sum_{p,q,r=1}^{d}(\sigma_{pq}(s,v)-%
\sigma_{pq}(s,z))(m_{r}(s,v)-m_{r}(s,z))\qquad\;
\]
\begin{equation}
\times\frac{\partial^{3}\widetilde{p}(s,jh,v,z)}{\partial v_{p}\partial
v_{q}\partial v_{r}}+\frac{1}{2}\sum_{p,q,r,l=1}^{d}\sigma_{pq}(s,v)\frac{%
\partial\sigma_{rl}(s,v)}{\partial v_{p}}\frac{\partial^{3}\widetilde{p}%
(s,jh,v,z)}{\partial v_{q}\partial v_{r}\partial v_{l}}+(\leq2),
\label{MR-012a}
\end{equation}
where we denote by $(\leq2)$ the sum of terms containing the derivatives of $%
\widetilde{p}(s,jh,v,z)$ of the order less or equal than 2. Note that for a
constant $C$ $<\infty$ and any $0<\varepsilon<\frac{1}{2}$%
\[
\left|\frac{h}{2}(p\otimes_{h}H_{1})(0,jh,x,z)\right|\leq Ch\phi_{C,\sqrt{jh}%
}\left(z-x\right),\qquad\qquad\;\;\;
\]
\begin{equation}
\left|\frac{h}{2}(p\otimes_{h}A_{0})(0,jh,x,z)\right|\leq
C(\varepsilon)h^{1/2}j^{-(1/2-\varepsilon)}\phi_{C,\sqrt{jh}%
}\left(z-x\right).  \label{MR-012b}
\end{equation}
First inequality (\ref{MR-012b}) follows from (B1) and the well know
estimates for the diffusion density $p$ and for the kernel $H_{1}$ . The
second inequality (\ref{MR-012b}) follows from (B1), (\ref{MR-012a}) and the
following estimate
\[
\frac{h}{2}\sum_{i=0}^{j-1}h\left|\int_{R^{d}}p(0,ih,x,v)\frac{\partial^{3}%
\widetilde{p}(ih,jh,v,z)}{\partial v_{q}\partial v_{r}\partial v_{l}}%
dv\right|=\frac{h}{2}\sum_{i=0}^{j-1}h\left|\int_{R^{d}}\frac{\partial
p(0,ih,x,v)}{\partial v_{q}}\right.
\]
\[
\left.\frac{\partial^{2}\widetilde{p}(ih,jh,v,z)}{\partial v_{r}\partial
v_{l}}dv\right|\leq Ch^{1-\varepsilon}\sum_{i=0}^{j-1}h\frac{1}{\sqrt{ih}}%
\frac{1}{(jh-ih)^{1-\varepsilon}}\phi_{C,\sqrt{jh}}\left(z-x\right)\leq
Ch^{1/2}\;
\]
\begin{equation}
\times j^{-(1/2-\varepsilon)}B(\frac{1}{2},\varepsilon)\phi_{C,\sqrt{jh}%
}\left(z-x\right).\qquad\qquad\qquad\qquad\qquad\qquad\qquad\qquad\qquad%
\qquad  \label{MR-012c}
\end{equation}
Now we shall estimate the second summand in the right hand side of (\ref
{MR-010}). Clearly
\[
\lambda^{\prime\prime}(s)=\frac{\partial^{2}}{\partial s^{2}}%
p(0,s,x,v)H(s,jh,v,z)+2\frac{\partial}{\partial s}p(0,s,x,v)
\]
\begin{equation}
\times\frac{\partial}{\partial s}H(s,jh,v,z)+p(0,s,x,v)\frac{\partial^{2}}{%
\partial s^{2}}H(s,jh,v,z).\qquad\;\;  \label{MR-013}
\end{equation}
Using forward and backward Kolmogorov equations we get from (\ref{MR-013})
after long but simple calculations
\[
\frac{1}{2}\sum_{i=0}^{j-1}\int_{ih}^{(i+1)h}(u-ih)^{2}\int_{0}^{1}(1-%
\delta)\int_{R^{d}}\lambda^{\prime\prime}(s)\mid_{s=s_{i}}dvd\delta
du\qquad\qquad\qquad\qquad\qquad\;\;
\]
\begin{equation}
=\frac{1}{2}\sum_{i=0}^{j-1}\int_{ih}^{(i+1)h}(u-ih)^{2}\int_{0}^{1}(1-%
\delta)\sum_{k=1}^{4}\int_{R^{d}}p(0,s,x,v)A_{k}(s,jh,v,z)\mid_{s=s_{i}}dvd%
\delta du,  \label{MR-014}
\end{equation}
where

\[
A_{1}(s,jh,v,z)=(L^{3}-3L^{2}\widetilde{L}+3L\widetilde{L}^{2}-\widetilde{L}%
^{3})\widetilde{p}(s,jh,v,z),\;\;\;\;\,
\]
\[
A_{2}=(L_{1}H+2LH_{1})(s,jh,v,z),\qquad \qquad \qquad \qquad \qquad \,\,
\]
\[
A_{3}(s,jh,v,z)=[(L-\widetilde{L})\widetilde{L}_{1}+2(L_{1}-\widetilde{L}%
_{1})\widetilde{L}]\widetilde{p}(s,jh,v,z),
\]
\begin{equation}
A_{4}(s,jh,v,z)=H_{2}(s,jh,v,z).\qquad \qquad \qquad \qquad \qquad \;\;\;
\label{MR-014a}
\end{equation}
and
\[
H_{l}(s,t,v,z)=(L_{l}-\widetilde{L}_{l})\widetilde{p}(s,t,v,z)\;\;\;\qquad
\qquad \qquad \qquad \qquad \;
\]
\[
=\frac{1}{2}\sum_{i,j=1}^{d}\left( \frac{\partial ^{l}\sigma _{ij}(s,v)}{%
\partial s^{l}}-\frac{\partial ^{l}\sigma _{ij}(s,z)}{\partial s^{l}}\right)
\frac{\partial ^{2}\widetilde{p}(s,t,v,z)}{\partial v_{i}\partial v_{j}}%
\qquad \qquad
\]
\begin{equation}
+\sum_{i=1}^{d}\left( \frac{\partial ^{l}m_{i}(s,v)}{\partial s^{l}}-\frac{%
\partial ^{l}m_{i}(s,z)}{\partial s^{l}}\right) \frac{\partial \widetilde{p}%
(s,t,v,z)}{\partial v_{i}},\;\;\;l=1,2.\;  \label{MR-015}
\end{equation}
Using integration by parts and the definition (\ref{MR-014a}) of $%
A_{2},A_{3} $ and $A_{4}$ it is easy to get that for any $0<\varepsilon <1/2$
and for $k=2,3,4$%
\[
\frac{1}{2}\left|
\sum_{i=0}^{j-1}\int_{ih}^{(i+1)h}(u-ih)^{2}\int_{0}^{1}(1-\delta
)\int_{R^{d}}p(0,s,x,v)A_{k}(s,jh,v,z)\mid _{s=s_{i}}dvd\delta du\right|
\]
\begin{equation}
\leq C(\varepsilon )h^{3/2-\varepsilon }\phi _{C,\sqrt{jh}}\left( z-x\right)
.\qquad \qquad \qquad \qquad \qquad \qquad \qquad \qquad \qquad \qquad \;
\label{MR-015a}
\end{equation}
For $k=1$ we shall prove the following estimate for any $0<\varepsilon <%
\frac{1}{2}$%
\[
\frac{1}{2}\left|
\sum_{i=0}^{j-1}\int_{ih}^{(i+1)h}(u-ih)^{2}\int_{0}^{1}(1-\delta
)\int_{R^{d}}p(0,s,x,v)A_{1}(s,jh,v,z)\mid _{s=s_{i}}dvd\delta du\right|
\]
\begin{equation}
\leq C(\varepsilon )hj^{-(1/2-\varepsilon )}\phi _{C,\sqrt{jh}}\left(
z-x\right) .\qquad \qquad \qquad \qquad \qquad \qquad \qquad \qquad \qquad
\qquad \;  \label{MR-015b}
\end{equation}
Note that the function $A_{1}(s,jh,v,z)$ can be written as the following sum
\[
A_{1}(s,jh,v,z)=\frac{1}{8}\sum_{i,j,p,q,l,r=1}^{d}(\sigma _{ij}(s,v)-\sigma
_{ij}(s,z))(\sigma _{pq}(s,v)-\sigma _{pq}(s,z))(\sigma _{lr}(s,v)\qquad
\]
\[
-\sigma _{lr}(s,z))\frac{\partial ^{6}\widetilde{p}(s,jh,v,z)}{\partial
v_{i}\partial v_{j}\partial v_{p}\partial v_{q}\partial v_{l}\partial v_{r}}+%
\frac{3}{4}\sum_{i,j,p,q,l=1}^{d}(\sigma _{ij}(s,v)-\sigma
_{ij}(s,z))(\sigma _{pq}(s,v)\qquad \;\;\;
\]
\[
-\sigma _{pq}(s,z))(m_{l}(s,v)-m_{l}(s,z))\frac{\partial ^{5}\widetilde{p}%
(s,jh,v,z)}{\partial v_{i}\partial v_{j}\partial v_{p}\partial v_{q}\partial
v_{l}}+\frac{3}{4}\sum_{i,j,p,q,l,r=1}^{d}\sigma _{ij}(s,v)\frac{\partial
\sigma _{pq}(s,v)}{\partial v_{i}}
\]
\begin{equation}
(\sigma _{lr}(s,v)-\sigma _{lr}(s,z))\times \frac{\partial ^{5}\widetilde{p}%
(s,jh,v,z)}{\partial v_{j}\partial v_{p}\partial v_{q}\partial v_{l}\partial
v_{r}}+(\leq 4),\;\;\;\qquad \qquad \qquad \qquad \qquad \qquad \,
\label{MR-016}
\end{equation}
where we denote by $(\leq 4)$ the sum of terms containing the derivatives of
$\widetilde{p}(s,jh,v,z)$ of the order less or equal than 4. By (B1) and (%
\ref{MR-016}) it is clear that the estimate for the left hand side of (\ref
{MR-015a}) for $k=1$ will be the same (up to a constant) as for the
following sum for fixed $p,q,r,l$

\[
\frac{1}{2}\left|
\sum_{i=0}^{j-1}\int_{ih}^{(i+1)h}(u-ih)^{2}\int_{0}^{1}(1-\delta
)\int_{R^{d}}p(0,s,x,v)\frac{\partial ^{4}\widetilde{p}(s,jh,v,z)}{\partial
v_{p}\partial v_{q}\partial v_{l}\partial v_{r}}\mid _{s=s_{i}}dvd\delta
du\right|
\]
After integration by parts w.r.t. $v_{p}$ and with the substitution $%
hw=(u-ih)$ in each integral we obtain
\[
\frac{1}{2}\left|
\sum_{i=0}^{j-1}\int_{ih}^{(i+1)h}(u-ih)^{2}\int_{0}^{1}(1-\delta
)\int_{R^{d}}p(0,s,x,v)\frac{\partial ^{4}\widetilde{p}(s,jh,v,z)}{\partial
v_{p}\partial v_{q}\partial v_{l}\partial v_{r}}\mid _{s=s_{i}}dvd\delta
du\right| \qquad \;\;\;\;\;
\]
\[
=\frac{1}{2}\left|
\sum_{i=0}^{j-1}\int_{ih}^{(i+1)h}(u-ih)^{2}\int_{0}^{1}(1-\delta
)\int_{R^{d}}\frac{\partial p(0,s,x,v)}{\partial v_{p}}\frac{\partial ^{3}%
\widetilde{p}(s,jh,v,z)}{\partial v_{q}\partial v_{l}\partial v_{r}}\mid
_{s=s_{i}}dvd\delta du\right| \;\;\;\;
\]
\[
\leq Ch^{2}\phi _{C,\sqrt{jh}}\left( z-x\right)
\int_{0}^{1}w^{2}\int_{0}^{1}(1-\delta )\sum_{i=0}^{j-1}h\frac{1}{\sqrt{%
ih+\delta hw}}\frac{1}{[(j-i)h-\delta hw]^{3/2}}d\delta dw.\qquad \;\;
\]
\[
\leq Ch^{3/2-\varepsilon }\phi _{C,\sqrt{jh}}\left( z-x\right)
\int_{0}^{1}w^{2}\int_{0}^{1}(1-\delta )^{1/2-\varepsilon }\sum_{i=0}^{j-1}h%
\frac{1}{\sqrt{ih+\delta hw}}\frac{1}{[(j-\delta w)h-ih]^{1-\varepsilon }}%
d\delta dw
\]
\[
\leq Ch^{3/2-\varepsilon }\phi _{C,\sqrt{jh}}\left( z-x\right)
\int_{0}^{1}w^{2}dw\int_{0}^{1}(1-\delta )^{1/2-\varepsilon }d\delta
\int_{0}^{(j-1)h}\frac{dt}{\sqrt{t}[(j-1)h-t]^{1-\varepsilon }}\qquad \qquad
\;\;\;
\]
\begin{equation}
\leq Chj^{-(1/2-\varepsilon )}B(\frac{1}{2},\varepsilon )\phi _{C,\sqrt{jh}%
}\left( z-x\right) ,\qquad \qquad \qquad \qquad \qquad \qquad \qquad \qquad
\qquad \qquad \qquad \;\;\;\;  \label{MR-017}
\end{equation}
where $B(p,q)$ is a Beta function and $\phi _{C,\rho }\left( z-x\right) $ is
defined in Lemma \ref{Lemma-KM-2} . As we mentioned above (\ref{MR-015b})
follows now from (\ref{MR-017}). By (\ref{MR-010}), (\ref{MR-012}), (\ref
{MR-012b}), (\ref{MR-015a}) and (\ref{MR-015b}) we obtain for any $%
0<\varepsilon <\frac{1}{2}$ and $j=1,2,...n$%
\begin{equation}
\left| \left( p\otimes H-p\otimes _{h}H\right) (0,jh,x,z)\right| \leq
C(\varepsilon )h^{1/2}j^{-(1/2-\varepsilon )}\phi _{C,\sqrt{jh}}\left(
z-x\right)  \label{MR-018}
\end{equation}
We use now the following estimate for $\Phi (ih,i^{\prime }h,z,z^{\prime })$
proved in Konakov and Mammen (2002)
\begin{equation}
\left| \Phi (ih,i^{\prime }h,z,z^{\prime })\right| \leq C\frac{1}{\sqrt{%
i^{\prime }h-ih}}\phi _{C,\sqrt{i^{\prime }h-ih}}\left( z^{\prime }-z\right)
\label{MR-019}
\end{equation}
From (\ref{MR-009}), (\ref{MR-018}) and (\ref{MR-019}) we obtain
\begin{equation}
\left| (p-p^{d})\left( 0,nh,x,y\right) \right| \leq C(\varepsilon
)h^{1/2}n^{\varepsilon -1/2}\phi _{C,\sqrt{T}}\left( y-x\right) .
\label{MR-020}
\end{equation}
The last inequality proves (\ref{MR-005}) for $k=1.$The terms $S_{2}$ and $%
S_{3}$ can be handled in the same way as the terms $T_{2}$ and $T_{3}$ in
Konakov and Mammen (2000). This completes the proof of our main result.

\textbf{Remark 1.} In fact for $S_{1}$ we proved stronger result than the
estimate (\ref{MR-005}). We get the following representation
\[
(p-p^{d})\left(0,T,x,y\right)=\frac{h}{2}(p\otimes_{h}H_{1})(0,T,x,y)+\frac{h%
}{2}(p\otimes_{h}A_{0})(0,T,x,y)
\]
\[
+\frac{h}{2}(p\otimes_{h}H_{1}\otimes_{h}\Phi)\left(0,T,x,y\right)+\frac{h}{2%
}(p\otimes_{h}A_{0}\otimes_{h}\Phi)(0,T,x,y)\qquad\;
\]
\begin{equation}
+R(0,T,x,y),\qquad\qquad\qquad\qquad\qquad\qquad\qquad\qquad\qquad\qquad%
\qquad\qquad\;  \label{MR-021}
\end{equation}
where for any $0<\varepsilon<1/2$%
\[
\left|R(0,T,x,y)\right|\leq
C(\varepsilon)(h^{3/2-\varepsilon}+hn^{-(1/2-\varepsilon)}+h\sqrt{T})\phi_{C,%
\sqrt{T}}\left(y-x\right)
\]
\[
=C(\varepsilon)\phi_{C,\sqrt{T}}\left(y-x\right)o(h).\qquad\qquad\qquad%
\qquad\qquad\qquad\qquad\;\;\,
\]
This representation is useful to obtain a small time Edgeworth type
expansions for our model (see Konakov and Mammen (2005))

\textbf{Remark 2.} If $T=const.$ (without loss of generality we assume $T=1)$
than we can easily avoid the difficulties connected with singularity by
splitting the time interval {[}0,1{]} and by using an integration by parts.
For example
\[
\int_{0}^{1}du\int_{R^{d}}p(0,u,x,v)\frac{\partial ^{4}\widetilde{p}(u,1,v,z)%
}{\partial v_{p}\partial v_{q}\partial v_{l}\partial v_{r}}%
dv=\int_{0}^{1/2}...+\int_{1/2}^{1}...
\]
\[
=\int_{0}^{1/2}du\int_{R^{d}}p(0,u,x,v)\frac{\partial ^{4}\widetilde{p}%
(u,1,v,z)}{\partial v_{p}\partial v_{q}\partial v_{l}\partial v_{r}}dv
\]
\begin{equation}
+\int_{1/2}^{1}du\int_{R^{d}}\frac{\partial ^{4}p(0,u,x,v)}{\partial
v_{p}\partial v_{q}\partial v_{l}\partial v_{r}}\widetilde{p}(u,1,v,z)dv
\label{MR-022}
\end{equation}
and the derivatives in the right hand side of (\ref{MR-022}) are not
singular. The representation (\ref{MR-021}) remains true. All summands in
the right hand side of (\ref{MR-021}) are estimated from above in absolute
value by $Ch\phi _{C,1}\left( y-x\right) $ and for the remainder term $%
R(0,1,x,y)$ the following estimate holds
\begin{equation}
\left| R(0,1,x,y)\right| \leq Ch^{2}\phi _{C,1}\left( y-x\right) .
\label{MR-023}
\end{equation}

\bigskip

\bigskip \textbf{References}

Konakov, V., Mammen, E. (2000). Local limit theorems for transition
densities of Markov chains converging to diffusions. Probab. Theory Rel.
Fields. 117, 551-587.

Konakov, V., Mammen, E. (2002). Edgeworth type expansions for Euler schemes
for stochastic differential equations. Monte Carlo Methods and Applications.
8, 271-286.

Konakov, V., Mammen, E. (2005). Edgeworth type expansions for transition
densities of Markov chains converging to diffusions. Bernoulli, v.11, n. 4,
p.591 - 641.

Friedman, A. (1964). Partial differential equations of parabolic type.
Prentice-Hall, Englewood Cliffs, New Jersey.

\end{document}